\documentclass[draft,11pt,english]{article}%{report}

\usepackage{amsfonts,amsmath,amsthm,amscd,amssymb,latexsym,verbatim,
caption2}%cite,texdraw,floatflt,pb-diagram

\theoremstyle{plain}
\newtheorem{theorem}{Theorem}[section]
\newtheorem{lemma}{Lemma}[section]
\newtheorem{proposition}{Proposition}[section]
\newtheorem{corollary}{Corollary}[section]
\newtheorem{definition}{Definition}[section]
\theoremstyle{definition}

\newtheorem{remark}{Remark}[section]
\newcommand{\keywords}{\textbf{Key words. }\medskip}
\newcommand{\subjclass}{\textbf{MSC 2010. }\medskip}
\renewcommand{\abstract}{\textbf{Abstract. }\medskip}

\numberwithin{equation}{section}
\def\MYim{\mathop{\rm Im}\limits}

\begin{document}

\title{SOS Polynomials and Matrix Representations\\ of Rational Real Functions}

\author{M. F. Bessmertny\u{\i}}

%\shorttitle{Short paper title for the headers}

%\shortauthor{F. Author, S. Author}

%\date{}

\maketitle

\begin{abstract}
 The characteristic properties of Artin's denominators in Hilbert's 17th problem are obtained. It is proved that numerators of partial derivatives of a rational real function from the Nevanlinna class are SOS polynomials.
\end{abstract}
\medskip

\subjclass{32A08, 15A22, 47N70.}

\keywords{Hilbert's 17th Problem, sums of squares, Nevanlinna's functions.}

    %\chapter{Chapter Title}

    \section*{Introduction}\label{s:1}

    A polynomial $F(z)\in \mathbb{R}[z_{1},\ldots,z_{d}]$ is called positive semidefinite (PSD) if it takes only non-negative values on $\mathbb{R}^{d}$ and it is called a sum of squares (SOS) if there exist other polynomials $h_{j}$ so that $F=h_{1}^{2}+\cdots+h_{k}^{2}$. Every SOS polynomial is a PSD, but for $d\geq 2$ there exist a PSD not SOS polynomials. In general, the PSD polynomial is the sum of the squares of rational functions (Hilbert's 17th problem).

    By Artin's solution to Hilbert's $17$th problem (see, for exemple, \cite{uj19}), for a PSD polynomial $F(z)$ there exists $s(z)\in \mathbb{R}[z_{1},\ldots,z_{d}]$ such that $s(z)^{2}F(z)$ is a SOS polynomial. The polynomial $s(z)$ will be called \emph{a Artin denominator} of the PSD polynomial $F(z)$.  Artin's theorem says nothing about the properties of the polynomial $s(z)$. In section \ref{s:2} we will prove (Theorem \ref{the2.1}) that there exists \emph{a minimal Artin denominator} whose irreducible factors do not change sign on $\mathbb{R}^{d}$.

    A PSD polynomials often appear as numerators of partial derivatives of rational function $f(z)=p(z)/q(z)$ (partial Wronskians $W_{j}[q(z),p(z)]$) \cite{uj01,uj02,uj23}. If the rational function $f(z)=p(z)/q(z)$ with real coefficients satisfies the condition (functions of the Nevanlinna class)
      \begin{equation}\label{eq1.1}
      \MYim f(z)\geq0\quad\text{for}\quad \MYim z_{1}\geq 0,\;\ldots\;, \MYim z_{d}\geq 0,
      \end{equation}
    then
      \begin{equation}\label{eq1.2}
      F_{j}(z)=W_{j}[q(z),p(z)]=q(z)\frac{\partial p(z)}{\partial z_{j}}-p(z)\frac{\partial q(z)}{\partial z_{j}},\quad
      j=1,\ldots,d.
      \end{equation}
    are PSD polynomials.
    \medskip

    \noindent
    \textbf{Main Theorem}.
    \emph{If rational function $f(z)=p(z)/q(z)$ with real coefficients satisfying condition \emph{(\ref{eq1.1})} , then}
    \begin{equation}\label{eq1.3}
    W_{j}[q(z),p(z)]=q(z)\frac{\partial p(z)}{\partial z_{j}}-p(z)\frac{\partial q(z)}{\partial z_{j}},\quad j=1,\ldots,d
    \end{equation}
    \emph{are sums of squares of polynomials.}
    \medskip

    This paper is organized as follows. Section \ref{s:2} introduces the Artin minimum denominator and studies the properties of its irreducible factors.

    The study of PSD polynomials (\ref{eq1.2}) is based on Product Polarization Theorem (Theorem \ref{the3.1}). This implies a long-resolvent representation (Corollary \ref{cor3.1}) for a rational real function (see \cite{uj03,uj04}).

    The ambiguity of the matrix pencil in the product polarization theorem is studied in Section \ref{s:4}. The ambiguity lemma (Lemma \ref{lem4.1}) allows one to obtain new matrix representations from a given one. A long-resolvent representation with one non-negative definite matrix is obtained for a function with one non-negative partial derivative on $\mathbb{R}^{d}$ (Corollary \ref{cor4.1}).

    On the basis of the representations obtained, in Section \ref{s:5} we prove the main theorem.

    \section{Artin's Denominators of PSD Polynomial}\label{s:2}

   If $F(z)\in\mathbb{R}[\,z_{1},\ldots,z_{d}\,]$ is a SOS polynomial, then $s(z)^{2}F(z)$ is also a SOS polynomial for each polynomial $s(z)$. If $F(z)$ is not representable as a sum of squares of polynomials, then the question arises: \emph {for which $s(z)$ is the polynomial $s(z)^{2}F(z)$ also not SOS polynomial?}

   \begin{proposition}\label{p2.1}
   \emph {(\cite{uj15}, Lemma 2.1).} Let $F(x)$ be a PSD not SOS polynomial and $s(x)$ an irreducible indefinite polynomial in $\mathbb{R}[\,x_{1},\ldots,x_{d}\,]$.
   Then $s^{2}F$ is also a PSD not SOS polynomial.
   \end{proposition}

  \begin{proof}
  Clearly $s^{2}F$ is PSD. If $s^{2}F=\sum_{k} h_{k}^{2}$, then for every real tuple $a$ with $s(a)=0$, it follows that $s^{2}F(a)=0$.
  This implies $h_{k}(a)^{2}=0$ \quad $\forall k$. So on the real manifold $s=0$, we have $h_{k}=0$ as well. So (see \cite{uj09}, Theorem 4.5.1)
  for each $k$, there exists $g_{k}$ so that $h_{k}=sg_{k}$. This gives $F=\sum_{k} g_{k}^{2}$, which is a contradiction.	
  \end{proof}

   \begin{proposition}\label{p2.2}
   Suppose $s^{2}F$ is a SOS polynomial. If each irreducible factor of $s$ is an indefinite polynomial, then $F$ is also a SOS polynomial.
   \end{proposition}

  \begin{proof}
  Suppose $F$ is a PSD not SOS polynomial. Let $s=s_{1}\cdots s_{m}$ be the decomposition of $s$ into irreducible factors.
  Successively applying Proposition \ref{p2.1} to the polynomials
    $$
    F_{1}=s_{1}^{2}F,\; F_{2}=s_{2}^{2}F_{1},\,\ldots\,,F_{m}=s_{m}^{2}F_{m-1},
    $$
  we obtain $F_{m}=s^{2}F$ is a PSD not SOS polynomial. Contradiction.
  \end{proof}

  \begin{definition}\label{d2.1}
  \emph{ Artin's denominator $s$ of a PSD polynomial $F$ is called \emph{an Artin minimal denominator}, if $\widehat{s}=s/s_{j}$ is
  not Artin's denominator of $F$ for each irreducible factor $s_{j}$ of $s$.}
  \end{definition}

   \begin{theorem}\label{the2.1}
   Each PSD not SOS polynomial $F(z)$ has a non-constant Artin minimum denominator $s(z)$. The irreducible factors of Artin's minimum
   denominator do not change sign on $\mathbb{R}^{d}$.
   \end{theorem}

   \begin{proof}
   By Artin's theorem, there exists $r(z)\in \mathbb{R}[z_{1},\ldots,z_{d}]$ for which $r(z)^{2}F(z)$ is a SOS. Each irreducible factor of
   $r(z)$ is either indefinite or does not change sign on $\mathbb{R}^{d}$. Then $r(z)=r_{0}(z)s(z)$, where all irreducible factors of the polynomial
   $s(z)$ do not change sign on $\mathbb{R}^{d}$, and irreducible factors of  $r_{0}(z)$ are indefinite. Consider the polynomial
   $F_{1}(z)=s(z)^{2}F(z)$. By condition, $r_{0}^{2}F_{1}=r^{2}F$ is a SOS. Since every irreducible factor of $r_{0}(z)$ is indefinite, we see that
   $F_{1}$ is a SOS (Proposition \ref{p2.2}). Then $s(z)$ is also Artin's denominator for $F$.
   Let $s_{0}$ be some irreducible factor of $s$. If $s/s_{0}$ remains the Artin denominator for  $F$, then the factor $s_{0}$ is
   removed from $s$. Removing all ``excess" irreducible factors from $s$, we obtain Artin's denominator with the required properties.
   \end{proof}

   To conclude this section, consider some properties of the irreducible factors of Artin's minimum denominator.

   \begin{proposition}\label{pro2.3}
   Let $d\geq 2$ and $s(z)\in\mathbb{R}[z_{1},\ldots,z_{d}]$ be an irreducible polynomial that does not change sign on $\mathbb{R}^{d}$. If $\partial s/\partial z_{d}\neq 0$, then there exist $\widehat{x}\in \mathbb{R}^{d-1}$ and $z_{d}\in\mathbb{C},\;\MYim z_{d}>0$ such that $s(\widehat{x},z_{d})=0$.
   \end{proposition}

   \begin{proof}
   Suppose for each fixed $\widehat{x}\in \mathbb{R}^{d-1}$ the polynomial $s(\widehat{x},z_{d})$ has only real zeros. Then the equation $s=0$ defines a real manifold of dimension $d-1$. By Theorem 4.5.1 \cite{uj09}, the ideal generated by the polynomial $s$ is a real, and $s(z)$ is indefinite. Contradiction.

   The complex zeros of $s(\widehat{x},z_{d})$  form complex conjugate pairs. Then for some $\widehat{x}\in\mathbb{R}^{d-1}$ there exist $z_{d}\in\mathbb{C},\;\MYim z_{d}>0$ such that $s(\widehat{x},z_{d})=0$.
   \end{proof}

   \begin{corollary}\label{cor2.1}
   If $s(z)\in\mathbb{R}[z_{1},\ldots,z_{d}]$ is a non-constant irreducible polynomial does not change sign on $\mathbb{R}^{d}$, then there exists $z'=(z'_{1},\ldots,z'_{d})$, $\MYim z'_{j}>0$, $j=1,\ldots,d$, such that $s(z')=0$.
   \end{corollary}

   \begin{proof}
   Without loss of generality, we can assume $\partial s/\partial z_{d}\neq 0$. By Proposition \ref{pro2.3}, $s(\widehat{x},z_{d})=0$ for $\widehat{x}=(x_{1},\ldots,x_{d-1})\in\mathbb{R}^{d-1}$, $\MYim z_{d}>0$.
   If $y_{1}>0,\ldots,y_{d-1}>0$ are small enough, then there is $z'_{d}$, $\MYim z'_{d}>0$ for which $s(x_{1}+iy_{1},\ldots,x_{d-1}+iy_{d-1},z'_{d})=0$.
   \end{proof}

    \section{Representation of Product of Polynomials}\label{s:3}

    Let $\delta_{j}$ be a non-negative integers. If $\alpha=(\delta_{1},\ldots,\delta_{d})$ and  $z=(z_{1},\ldots,z_{d})$, then by $z^{\alpha}$ denote a monomial $z_{1}^{\delta_{1}}\cdots z_{d}^{\delta_{d}}$.

    \begin{theorem}\label{the3.1} \emph{(\textbf{Product Polarization Theorem})}.
    Let $q(z),\,p(z)$, where $z\in\mathbb{C}^{d}$, be polynomials with real coefficients, and let $\Psi(z)=(z^{\alpha_{1}},\ldots,z^{\alpha_{N}})$ be the row vector of all monomials satisfying the conditions
    \begin{equation}\label{eq3.1}
      \begin{array}{c}
      \quad\quad\;\deg z^{\alpha_{j}}\leq\max\{\deg q(z),\,\deg p(z)\},\qquad\qquad\qquad\quad\\
      \deg_{z_{k}}z^{\alpha_{j}}\leq\max\{\deg_{z_{k}} q(z),\,\deg_{z_{k}} p(z)\},\;k=1,\,\ldots\,,d;\quad
      \end{array}
    \end{equation}
    then \emph{(i)} there exist real symmetric $N\times N$-matrices $B_{j}$, $j=0,1,\,\ldots\,,d$ such that
        \begin{equation}\label{eq3.2} q(\zeta)p(z)=\Psi(\zeta)(B_{0}+z_{1}B_{1}+\cdots+z_{d}B_{d})\Psi(z)^{T},\quad \forall\zeta,z\in \mathbb{C}^{d},
        \end{equation}
    \emph{(ii)} if $q(z)\partial p(z)/\partial z_{d}-p(z)\partial q(z)/\partial z_{d}$ is the sum of squares of  polynomials, then the matrix $B_{d}$ can be chosen positive semidefinite.
    \end{theorem}
    In this section we prove assertion (i) of Theorem \ref{the3.1}. The proof of (ii) will be given at the end of section \ref{s:4}.
    \medskip

    We need a few lemmas.
    A monomial $z^{\alpha}$ is called  \emph{multi-affine} if $\deg_{z_{k}}z^{\alpha}\leq 1$ for every $k=1,\ldots, d$.
    \begin{lemma}\label{lem3.1}
    Let $k\geq0$ be an integer. If $\zeta^{\mu_{1}}=\zeta_{2}\zeta_{4}\cdots\zeta_{2k}$ and $\zeta^{\nu}=\zeta_{1}\zeta_{3}\cdots\zeta_{2k+1}$, then there exist real symmetric $(2k+1)\times(2k+1)$ matrices $C_{j}$, $j=1,2,\,\ldots\,,2k+1$  and multiaffine monomials $\{\zeta^{\mu_{j}}\}_{j=2}^{2k+1}$
    of degree $k$ in variables $\zeta_{1},\zeta_{1},\,\ldots\,,\zeta_{2k+1}$ such that
       \begin{equation}\label{eq3.3}
       (\zeta_{1}C_{1}+\zeta_{2}C_{2}+\cdots+\zeta_{2k+1}C_{2k+1})
         \begin{pmatrix}
         \zeta^{\mu_{1}}\\
         \zeta^{\mu_{2}}\\
         \vdots\\
         \zeta^{\mu_{2k+1}}
         \end{pmatrix}=
           \begin{pmatrix}
           \zeta^{\nu}\\
           0\\
           \vdots\\
           0
           \end{pmatrix}.
       \end{equation}
    \end{lemma}

    \begin{proof}
    If $k=0$, then $\zeta^{\mu_{1}}=1$ (empty product) and $\zeta^{\nu}=\zeta_{1}$. We have $\zeta^{\nu}=\zeta_{1}\zeta^{\mu_{1}}$, and the matrix pencil has the size $1\times 1$.
    For $k\geq 1$ multiaffine monomials $\{\zeta^{\mu_{j}}\}_{j=2}^{2k+1}$ are defined by the relations
      $$
      \zeta^{\mu_{2}}=\zeta_{3}\zeta_{5}\cdots\zeta_{2k+1},\quad
      \zeta^{\mu_{j}}=\zeta_{j-2}\zeta^{\mu_{j-2}}/\zeta_{j-1},\;
      j=3,4,\,\ldots\,,2k+1.
      $$
    Note that $\zeta^{\mu_{2k+1}}=\zeta^{\nu}/\zeta_{2k+1}$, $\zeta^{\mu_{2k}}=\zeta_{2k+1}\zeta^{\mu_{1}}/\zeta_{2k}$. Let us define the matrix pencil $C(\zeta)=\{c_{ij}(\zeta)\}_{i,j=1}^{2k+1}$:
      $$
      c_{ij}(\zeta)=
        \begin{cases}
        (-1)^{\max\{i,j\}}\zeta_{\min\{i,j\}}/2,& \text{if $|i-j|=1$},\\
        \zeta_{\max\{i,j\}}/2,& \text{if $|i-j|=2k$},\\
        0,& \text{otherwise}.
        \end{cases}
      $$
    It is easy to see $\overline{C(\overline{\zeta})}=C(\zeta)=C(\zeta)^{T}$. Let us calculate the components $b_{i}=b_{i}(\zeta)$ of the right-hand side of (\ref{eq3.3}).
      $$
      b_{1}(\zeta)=\sum_{j=1}^{2k+1}c_{1j}(\zeta)\zeta^{\mu_{j}}=c_{12}(\zeta)\zeta^{\mu_{2}}+
      c_{1,2k+1}(\zeta)\zeta^{\mu_{2k+1}}=\zeta^{\nu}.
      $$
    For $2\leq i\leq 2k$ we obtain
      \begin{multline}\nonumber
      b_{i}(\zeta)=\sum_{j=1}^{2k+1}c_{ij}(\zeta)\zeta^{\mu_{j}}=c_{i,i-1}(\zeta)\zeta^{\mu_{i-1}}+
      c_{i,i+1}(\zeta)\zeta^{\mu_{i+1}}=\\
      (-1)^{i}(\zeta_{i-1}\zeta^{\mu_{i-1}}-\zeta_{i}\zeta^{\mu_{i+1}})/2=
      (-1)^{i}(\zeta_{i-1}\zeta^{\mu_{i-1}}-\zeta_{i-1}\zeta^{\mu_{i-1}})/2=0.
      \end{multline}
    For $i=2k+1$ we get
      \begin{multline}\nonumber
      b_{2k+1}(\zeta)=\sum_{j=1}^{2k+1}c_{2k+1,j}(\zeta)\zeta^{\mu_{j}}=
      c_{2k+1,1}(\zeta)\zeta^{\mu_{1}}+c_{2k+1,2k}(\zeta)\zeta^{\mu_{2k}}=\\
      (\zeta_{2k+1}\zeta^{\mu_{1}}-\zeta_{2k}\zeta^{\mu_{2k}})/2=
      (\zeta_{2k+1}\zeta^{\mu_{1}}-\zeta_{2k}\zeta_{2k+1}\zeta^{\mu_{1}}/\zeta_{2k})/2=0,\quad
      \end{multline}
    what was required.
    \end{proof}

    \begin{lemma}\label{lem3.2}
    Let $z^{\alpha_{1}}$, $z^{\beta}$ be arbitrary monomials in variables $z_{1},\,\ldots\,,z_{d}$. Then there are real symmetric matrices $D_{j}$, $j=0,1,\,\ldots\,,d$ such that
       \begin{equation}\label{eq3.4}
       (D_{0}+z_{1}D_{1}+\cdots+z_{d}D_{d})
         \begin{pmatrix}
         z^{\alpha_{1}}\\
         z^{\alpha_{2}}\\
         \vdots\\
         z^{\alpha_{N}}
         \end{pmatrix}=
           \begin{pmatrix}
           z^{\beta}\\
           0\\
           \vdots\\
           0
           \end{pmatrix},
       \end{equation}
    where $\{z^{\alpha_{j}}\}_{j=2}^{N}$ are all monomials satisfying the conditions
       \begin{equation}\label{eq3.5}
       \begin{array}{c}
       \quad\quad\quad\deg z^{\alpha_{j}}\leq\max \{\deg z^{\alpha_{1}},\,\deg z^{\beta}\},\qquad\qquad\qquad\qquad\\
       \deg_{z_{k}} z^{\alpha_{j}}\leq\max \{\deg_{z_{k}} z^{\alpha_{1}},\,\deg_{z_{k}}z^{\beta}\},\;
       k=1,\,\ldots\,,d.\quad\;\,
       \end{array}
       \end{equation}
    \end{lemma}

    \begin{proof}
    Let $l,m\in \mathbb{N}$ be such that
       $$
       \deg z_{0}^{l}z^{\alpha_{1}}=n,\quad \deg z_{0}^{m}z^{\beta}=n+1,
       $$
    where $z_{0}$ is new variable.
    Let $z^{\gamma}$ be the greatest common divisor of the monomials $z_{0}^{l}z^{\alpha_{1}}$ and $z_{0}^{m}z^{\beta}$. Then the monomials $z_{0}^{l}z^{\alpha_{1}}/z^{\gamma}$ and $z_{0}^{m}z^{\beta}/z^{\gamma}$ do not contain the same variables and have degrees $k$ and $k+1$.

    \noindent
    By Lemma \ref{lem3.1}, there exist real symmetric matrices $C_{j}$, $j=1,\ldots,2k+1$ such that
       \begin{equation}\label{eq3.6}
       (\zeta_{1}C_{1}+\cdots+\zeta_{2k+1}C_{2k+1})
         \begin{pmatrix}
         \zeta_{2}\zeta_{4}\cdots\zeta_{2k}z^{\gamma}\\
         \zeta^{\mu_{2}}z^{\gamma}\\
         \vdots\\
         \zeta^{\mu_{2k+1}}z^{\gamma}
         \end{pmatrix}=
           \begin{pmatrix}
           \zeta_{1}\zeta_{3}\cdots\zeta_{2k+1}z^{\gamma}\\
           0\\
           \vdots\\
           0
           \end{pmatrix},
       \end{equation}
    where $\{\zeta^{\mu_{j}}\}_{j=2}^{2k+1}$ are monomials of degree $k$ in variables $\zeta_{1},\zeta_{2},\,\ldots\,,\zeta_{2k+1}$.

    In (\ref{eq3.6}) we replace the variables $\zeta_{2}\zeta_{4}\cdots\zeta_{2k}$ by the variables of the monomial $z_{0}^{l}z^{\alpha_{1}}/z^{\gamma}$, and the variables $\zeta_{1}\zeta_{3}\cdots\zeta_{2k+1}$ by the variables of monomial $z_{0}^{m}z^{\beta}/z^{\gamma}$, so that
      $$
      \zeta_{2}\zeta_{4}\cdots\zeta_{2k}z^{\gamma}\mapsto z_{0}^{l}z^{\alpha_{1}},\quad
       \zeta_{1}\zeta_{3}\cdots\zeta_{2k+1}z^{\gamma}\mapsto z_{0}^{m}z^{\beta}.
      $$
    Setting $z_{0}=1$, we obtain
      \begin{equation}\label{eq3.7}
      (D'_{j_{0}}+z_{j_{1}}D'_{j_{1}}+\cdots+z_{j_{r}}D'_{j_{r}})
          \begin{pmatrix}
          z^{\alpha_{1}}\\
          z^{\widehat{\alpha}_{2}}\\
          \vdots\\
          z^{\widehat{\alpha}_{2k+1}}
          \end{pmatrix}=
             \begin{pmatrix}
             z^{\beta}\\
             0\\
             \vdots\\
             0
             \end{pmatrix},
      \end{equation}
    where $z_{j_{1}},\,\ldots\,,z_{j_{r}}$ are the variables of monomials $z_{0}^{l}z^{\alpha_{1}}/z^{\gamma}$, $z_{0}^{m}z^{\beta}/z^{\gamma}$ different from $z_{0}$. The matrices $D'_{j_{k}}$ are sums of the corresponding matrices $C_{j}$ from (\ref{eq3.6}). It is easy to see (\ref{eq3.5}) holds.

    Let $\{z^{\alpha_{j}}\}_{j=1}^{N}$, $z^{\alpha_{i}}\neq z^{\alpha_{j}}$, $i\neq j$, be all monomials satisfying (\ref{eq3.5}). Then there exists a $(2k+1)\times N$ matrix $B$ such that
       $$
       \begin{pmatrix}
       z^{\alpha_{1}}\\
       z^{\widehat{\alpha}_{2}}\\
       \vdots\\
       z^{\widehat{\alpha}_{2k+1}}
       \end{pmatrix}=
          \begin{pmatrix}
          1          &   0        & \cdots & 0\\
          b_{21}     & b_{22}     & \cdots & b_{2N}\\
          \vdots     & \vdots     & \vdots & \vdots \\
          b_{2k+1,1} & b_{2k+1,1} & \cdots & b_{2k+1,N}
          \end{pmatrix}
            \begin{pmatrix}
            z^{\alpha_{1}}\\
            z^{\alpha_{2}}\\
            \vdots\\
            z^{\alpha_{N}}
            \end{pmatrix}.
       $$
    From (\ref{eq3.7}), introducing the notation
       $$
       B^{T}(D'_{j_{0}}+z_{j_{1}}D'_{j_{1}}+\cdots+z_{j_{r}}D'_{j_{r}})B=D_{0}+z_{1}D_{1}+\cdots+z_{d}D_{d},
       $$
    where part of the matrices $D_{k}$ is equal to zero, we obtain (\ref{eq3.4}).
    \end{proof}

    \noindent
    \emph{Proof of assertion \emph{(i)} of Theorem \ref{the3.1}}. Suppose
      $$
      q(z)=\sum_{j=1}^{N}a_{j}z^{\alpha_{j}},\quad p(z)=\sum_{\nu=1}^{M}b_{\nu}z^{\beta_{\nu}},\quad a_{j},b_{\nu}\in \mathbb{R}.
      $$
    By Lemma \ref{lem3.2},
    for fixed monomials $z^{\alpha_{j}}$ and $z^{\beta_{\nu}}$ there exists a symmetric real matrix pencil $D_{j\nu}(z)$ such that
      \begin{equation}\nonumber
      D_{j\nu}(z)
        \begin{pmatrix}
        z^{\alpha_{1}}\\
        \vdots\\
        z^{\alpha_{j}}\\
        \vdots\\
        z^{\alpha_{N}}
        \end{pmatrix}=
          \begin{pmatrix}
          0\\
          \vdots\\
          z^{\beta_{\nu}}\\
          \vdots\\
          0
          \end{pmatrix} - j\text{-th row},\quad j=1,\ldots,N,\;\nu=1,\ldots,M,
      \end{equation}
    where $\{z^{\alpha_{j}}\}_{j=1}^{N}$ are all monomials satisfying the conditions
       \begin{equation}\nonumber
       \begin{array}{c}
       \deg z^{\alpha_{j}}\leq\max\{\deg q(z),\,\deg p(z)\},\qquad\quad\;\\
    \deg_{z_{k}}z^{\alpha_{j}}\leq\max\{\deg_{z_{k}}q(z),\,\deg_{z_{k}}p(z)\},\;  k=1,\ldots,d.
   \end{array}
   \end{equation}
    We define
    $$
    B(z)=B_{0}+z_{1}B_{1}+\cdots+z_{d}B_{d}=\sum_{j=1}^{N}a_{j}\sum_{\nu=1}^{M}b_{\nu}D_{j\nu}(z).
    $$
    It's easy to see
      $$
      B(z)
        \begin{pmatrix}
        z^{\alpha_{1}}\\
        z^{\alpha_{2}}\\
        \vdots\\
        z^{\alpha_{N}}
        \end{pmatrix}=
          \begin{pmatrix}
          a_{1}p(z)\\
          a_{2}p(z)\\
          \vdots\\
          a_{N}p(z)
          \end{pmatrix}.
      $$
    Multiplying the last identity on the left by the row vector $\Psi(\zeta)=(\zeta^{\alpha_{1}},\ldots,\zeta^{\alpha_{N}})$, we obtain (\ref{eq3.2}).
    \qed

    \begin{corollary}\label{cor3.1} \emph{\textbf{(Long-resolvent Representation)}}.
    Let $z\in \mathbb{C}^{d}$.
    For a rational function $f(z)$ with real coefficients there exists $M\times M$-matrix pencil
      \begin{equation}\label{eq3.8}
         A_{0}+z_{1}A_{1}+\cdots+z_{d}A_{d}=
       \begin{pmatrix}
       A_{11}(z) & A_{12}(z)\\
       A_{21}(z) & A_{22}(z)
       \end{pmatrix},
      \end{equation}
    where  $\overline{A}_{j}=A_{j}=A_{j}^{T}$, $j=0,1,\ldots,d$, such that
      \begin{equation}\label{eq3.9}
      f(z)=A_{11}(z)-A_{12}(z)A_{22}(z)^{-1}A_{21}(z).
      \end{equation}
    \end{corollary}

    \begin{proof}
    Suppose $f(z)=p(z)/q(z)$, where $q(z)=\sum_{j=1}^{N}a_{j}z^{\alpha_{j}}$. Without loss of generality, we can assume $a_{1}\neq0$. By Theorem \ref{the3.1} (i), there exists a real symmetric $N\times N$-matrix pencil $B(z)=B_{0}+z_{1}B_{1}+\cdots+z_{d}B_{d}$ such that
       \begin{equation}\label{eq3.10}
       B(z)
         \begin{pmatrix}
         z^{\alpha_{1}}\\
         z^{\alpha_{2}}\\
         \vdots\\
         z^{\alpha_{N}}
         \end{pmatrix}=
           \begin{pmatrix}
           a_{1}p(z)\\
           a_{2}p(z)\\
           \vdots\\
           a_{N}p(z)
           \end{pmatrix}.
       \end{equation}
    If
      $$
      Q=
       \begin{pmatrix}
       a_{1} & a_{2}  & \cdots & a_{N}  \\
        0    &  1     & \cdots &   0    \\
      \vdots & \vdots & \ddots & \vdots \\
        0    &  0     & \cdots &   1
       \end{pmatrix},
      $$
    then assuming $ \widetilde{A}(z)=Q^{-1T}B(z)Q^{-1}$, we get
      \begin{equation}\label{eq3.11}
      \begin{pmatrix}
               A_{11}(z)        & \widetilde{A}_{12}(z)\\
      \widetilde{A}_{12}(z)^{T} & \widetilde{A}_{22}(z)
      \end{pmatrix}
         \begin{pmatrix}
             q(z)     \\
         z^{\alpha_{2}}\\
         \vdots\\
         z^{\alpha_{N}}
         \end{pmatrix}=
           \begin{pmatrix}
             p(z) \\
              0   \\
           \vdots \\
              0
           \end{pmatrix}.
      \end{equation}
    Suppose $\det \widetilde{A}_{22}(z)=0$. Let $A_{22}(z)$ be the principal submatrix of the maximum size of matrix $\widetilde{A}_{22}(z)$ such that $\det A_{22}(z)\neq0$. Without loss of generality, we can assume that $A_{22}(z)$ is located in the upper left corner of the matrix $\widetilde{A}_{22}(z)$. (\ref{eq3.11}) rewrite in the form
      $$
      \widetilde{A}(z)
      \begin{pmatrix}
            q(z)    \\
        \psi_{2}(z) \\
        \psi_{3}(z)
        \end{pmatrix}=
      \begin{pmatrix}
          A_{11}(z) &     A_{12}(z)  & A_{13}(z) \\
      A_{12}(z)^{T} &     A_{22}(z)  & A_{23}(z) \\
      A_{13}(z)^{T} &  A_{13}(z)^{T} & A_{33}(z)
      \end{pmatrix}
        \begin{pmatrix}
            q(z)    \\
        \psi_{2}(z) \\
        \psi_{3}(z)
        \end{pmatrix}=
          \begin{pmatrix}
          p(z) \\
            0  \\
            0
          \end{pmatrix}.
      $$
    Let's put
       $$
       R(z)=
         \begin{pmatrix}
         1 & 0 &           0            \\
         0 & I & A_{22}(z)^{-1}A_{23}(z)\\
         0 & 0 &           I
         \end{pmatrix}.
       $$
    Assuming $\widehat{A}(z)=R(z)^{-1T}\widetilde{A}(z)R(z)^{-1}$, from the last relation we obtain
       \begin{equation}\label{eq3.12}
       \begin{pmatrix}
            A_{11}(z)     &  A_{12}(z) & C_{12}(z)\\
            A_{12}(z)^{T} &  A_{22}(z) &     0    \\
          C_{12}(z)^{T}   &      0     &     0
       \end{pmatrix}
          \begin{pmatrix}
              q(z)     \\
            \varphi(z) \\
            \psi_{3}(z)
          \end{pmatrix}=
            \begin{pmatrix}
              p(z) \\
                0  \\
                0
            \end{pmatrix},
       \end{equation}
    where $C_{12}(z)=A_{13}-A_{12}A_{22}^{-1}A_{23}$, $\varphi(z)=\psi_{2}+A_{22}^{-1}A_{23}\psi_{3}$.

    Since $q(z)\neq0$, then $C_{12}(z)^{T}=0$. Therefore, from (\ref{eq3.12}) we have
    $\varphi(z)=-A_{22}(z)^{-1}A_{12}(z)^{T}q(z)$, $p(z)=A_{11}(z)q(z)+A_{12}(z)\varphi(z)$, whence follows (\ref{eq3.9}).
    \end{proof}

    \begin{remark}\label{rem3.1}
    If one of the matrices $B_{k}$, $k=0,1,\ldots,d$ (see (\ref{eq3.10})) is positive semidefinite, then the corresponding matrix in the  matrix pencil $A(z)$ (\ref{eq3.8}) will also be positive semidefinite.
    \end{remark}

     To prove assertion (ii) of Theorem \ref{the3.1}, we need to describe the ambiguities of representation (\ref{eq3.2}).

    \section{Ambiguities of the matrix pencil}\label{s:4}

    If $p(z),q(z)$, $z\in\mathbb{C}^{d}$, are polynomials with real coefficients and
      \begin{equation}\label{eq4.1}
      \max \{\deg p,\,\deg q\}=n_{0};\quad\max \{\deg_{z_{k}} p,\,\deg_{z_{k}}q\}=n_{k}
      \end{equation}
    for $k=1\ldots,d$, then
    by $\Psi(z)=(z^{\alpha_{1}},\ldots,z^{\alpha_{N}})$  we denote the row vector of all monomials $z^{\alpha_{j}}=z_{1}^{\delta_{1}}\cdots z_{d}^{\delta_{d}}$ satisfying the conditions
      \begin{equation}\label{eq4.2}
      \deg z^{\alpha_{j}}\leq n_{0},\quad \deg z^{\alpha_{j}}\leq n_{k},\quad k=1\ldots,d.
      \end{equation}

    \begin{proposition}\label{pro4.2}
    Let $p(z),q(z)$ be polynomials satisfying \emph{(\ref{eq4.1})} and $A(z)=A_{0}+z_{1}A_{1}+\cdots+z_{d}A_{d}$ be a symmetric matrix pencil such that
      \begin{equation}\label{eq4.3}
      q(\zeta)p(z)=\Psi(\zeta)(A_{0}+z_{1}A_{1}+\cdots+z_{d}A_{d})\Psi(z)^{T},\quad \zeta,z\in\mathbb{C}^{d},
      \end{equation}
    where $\Psi(z)=(z^{\alpha_{1}},\ldots,z^{\alpha_{N}})$ is the row vector of all monomials $z^{\alpha_{j}}$ satisfying the conditions \emph{(\ref{eq4.2})}. Then for every $k=1,\ldots,d$
    \begin{multline}\nonumber
    \emph{(i)}.\quad
      \Psi(z)A_{k}\Psi(z)^{T}=q(z)\frac{\partial p}{\partial z_{k}}(z)-p(z)\frac{\partial q}{\partial z_{k}}(z)=W_{k}[q(z),p(z)].
    \end{multline}
    \begin{multline}\nonumber
    \emph{(ii)}.\quad A_{k}\frac{\partial^{n_{k}} \Psi^{T}(z)}{\partial z^{n_{k}}}=0.\\
    \end{multline}

    \end{proposition}
    \begin{proof}
    (i). From (\ref{eq4.3}) we obtain
      \begin{multline}\nonumber
      q(\zeta)\frac{\partial p}{\partial z_{k}}(z)-p(z)\frac{\partial q}{\partial \zeta_{k}}(\zeta)=\\
      \Psi(\zeta)A_{k}\Psi(z)^{T}+\Psi(\zeta)A(z)\frac{\partial\Psi}{\partial z_{k}}(z)-
      \frac{\partial\Psi}{\partial\zeta_{k}}(\zeta)A(z)\Psi(z)^{T}.
      \end{multline}
    Since $A(z)^{T}=A(z)$, then assuming $\zeta=z$ we get
      $$
      q(z)\frac{\partial p}{\partial z_{k}}(z)-p(z)\frac{\partial q}{\partial z_{k}}(z)=
      \Psi(z)A_{k}\Psi(z)^{T},\quad k=1,\ldots,d.
      $$

    \noindent
    (ii). Differentiating (\ref{eq4.3}) $n_{k}+1$ times with respect to  $z_{k}$, we obtain
      $$
      q(\zeta)\frac{\partial^{n_{k}+1}p(z)}{\partial z_{k}^{n_{k}+1}}=
      (n_{k}+1)\Psi(\zeta)A_{k}\frac{\partial^{n_{k}}\Psi(z)^{T}}{\partial z^{n_{k}}}+
      \Psi(\zeta)A(z)\frac{\partial^{n_{k}+1}\Psi(z)^{T}}{\partial z^{n_{k}+1}}.
      $$
    Since
      $$
      \frac{\partial^{n_{k}+1}p(z)}{\partial z_{k}^{n_{k}+1}}=0\quad\text{and}\quad
      \frac{\partial^{n_{k}+1}\Psi(z)^{T}}{\partial z^{n_{k}+1}}=0,
      $$
    then
      $$
      A_{k}\frac{\partial^{n_{k}}\Psi(z)^{T}}{\partial z^{n_{k}}}=0,\quad\text{where}\quad
      \frac{\partial^{n_{k}}\Psi(z)^{T}}{\partial z^{n_{k}}}\neq 0,
      $$
    what was required.
    \end{proof}

    \begin{remark}\label{rem4.1}
    If there are two representations (\ref{eq4.3}) with pencils $A(z)$ and $B(z)$, then $S(z)=A(z)-B(z)$ satisfies the conditions
      $$
      (S_{0}+z_{1}S_{1}+\cdots+z_{d}S_{d})\Psi(z)^{T}=0,
      $$
       $$
       \Psi(z)S_{k}\Psi(z)^{T}=0,\quad S_{k}\frac{\partial^{n_{k}}\Psi(z)^{T}}{\partial z_{k}^{n_{k}}}=0,\quad k=0,1,\ldots,d.
       $$
    \end{remark}

   \begin{lemma}\label{lem4.1} \emph{\textbf{(Ambiguity Lemma)}.}
   Let $S_{d}$ be a real symmetric $N\times N$-matrix for which
     \begin{equation}\label{eq4.4}
     \Psi(z)S_{d}\Psi(z)^{T}=0,\quad S_{d}\frac{\partial^{n_{d}}\Psi(z)^{T}}{\partial z_{d}^{n_{d}}}=0,
     \end{equation}
   where $\Psi(z)=(z^{\alpha_{1}},\,\ldots\,,z^{\alpha_{N}})$  is the row vector of all monomials $z^{\alpha_{j}}$ satisfying conditions \emph{(\ref{eq4.2})}.
   Then there exist real symmetric $N\times N$-matrices $S_{0},S_{1},\ldots,S_{d-1}$ such that
     \begin{equation}\label{eq4.5}
     (S_{0}+z_{1}S_{1}+\cdots+z_{d-1}S_{d-1}+z_{d}S_{d})\Psi(z)^{T}=0.
     \end{equation}
   \end{lemma}

   \begin{remark}\label{rem4.2}
   It is easy to see, the matrices $S_{k}$, $k=0,1,\ldots,d$, satisfy condition (\ref{eq4.5}) if and only if
     $$
     (z_{0}S_{0}+z_{1}S_{1}+\cdots+z_{d-1}S_{d-1}+z_{d}S_{d})\widetilde{\Psi}(z)^{T}=0,\quad z\in\mathbb{C}^{d+1},
     $$
   where $\widetilde{\Psi}(z)=(z^{\alpha_{1}},\ldots,z^{\alpha_{N}})$ is a row vector of all monomials $z^{\alpha_{j}}=z_{0}^{\delta_{0}}z_{1}^{\delta_{1}}\cdots z_{d}^{\delta_{d}}$ such that
     $$
     \deg z^{\alpha_{j}}=n_{0},\quad \deg_{z_{k}} z^{\alpha_{j}}\leq n_{k},\quad k=0,1,\ldots,d.
     $$
   In addition,\; $\Psi(z)S_{d}\Psi(z)^{T}=0$ \;iff\;  $\widetilde{\Psi}(z)S_{d}\widetilde{\Psi}(z)^{T}=0$. \qed
   \end{remark}

   \begin{remark}\label{rem4.3}
   Let $L$ be a linear space of real symmetric matrices $S=\{s_{ij}\}_{i,j=1}^{N}$ such that
      \begin{equation}\label{eq4.6}
      \widetilde{\Psi}(z)S\widetilde{\Psi}(z)^{T}=0.
      \end{equation}
   Since
      $$
      \widetilde{\Psi}(z)S\widetilde{\Psi}(z)^{T}=\sum_{i,j=1}^{N}s_{ij}z^{\alpha_{i}}z^{\alpha_{j}}
      =\sum_{k=1}^{M}c_{k}z^{\beta_{k}}=0,
      $$
   where $\beta_{k}\neq\beta_{\nu}$, $k\neq\nu$, then $L$ is the direct sum of subspaces $L_{\beta_{k}}$ consisting of matrices  $S_{\beta_{k}}=\{s_{ij}\}_{i,j=1}^{N}$ satisfying the conditions
     \begin{equation}\label{eq4.7}
     \sum_{\alpha_{i}+\alpha_{j}=\beta_{k}}s_{ij}=0\quad\text{and}\quad s_{ij}=0\quad\text{for}\quad
     \alpha_{i}+\alpha_{j}\neq\beta_{k},
     \end{equation}
   where $\alpha_{i}+\alpha_{j}$ is the component-wise sum of multi-indices.

   The subspace $L_{\beta_{k}}$ is non-zero if and only if there are at least $m\geq 2$ different non-zero elements $s_{ij}=s_{ji}$ such that $\alpha_{i}+\alpha_{j}=\beta_{k}$.
   In this case $\dim L_{\beta_{k}}=m-1$. \qed
   \end{remark}

    Let us find a basis in each subspace $L_{\beta_{k}}$.
    \begin{proposition}\label{pro4.2}
    Let $L$ be a linear space of real symmetric matrices satisfying condition \emph{(\ref{eq4.6})}, and let $L_{\beta_{k}}\neq\{0\}$ be a subspace of matrices $S_{\beta_{k}}=\{s_{ij}\}_{i,j=1}^{N}$ for which \emph{(\ref{eq4.7})} holds. In $L_{\beta_{k}}$ there exists a basis such that the nonzero elements of the basis matrices are at the intersection of rows and columns corresponding to monomials of a special form:
      \begin{equation}\label{eq4.8}
       (z_{r}^{2}z^{\gamma},\, z_{r}z_{l}z^{\gamma},\, z_{l}^{2}z^{\gamma})
         \begin{pmatrix}
         0  & 0 & -1\\
         0  & 2 &  0\\
         -1 & 0 &  0
         \end{pmatrix}
            \begin{pmatrix}
            z_{r}^{2}z^{\gamma}\\
            z_{r}z_{l}z^{\gamma}\\
            z_{l}^{2}z^{\gamma}
            \end{pmatrix},
      \end{equation}
        \begin{equation}\label{eq4.9}
         (z_{\mu}z^{\gamma_{1}},\,z_{\nu}z^{\gamma_{1}},\,z_{\nu}z^{\gamma_{2}},\, z_{\mu}z^{\gamma_{2}})
           \begin{pmatrix}
           0 &  0 & 1 &  0\\
           0 &  0 & 0 & -1\\
           1 &  0 & 0 &  0\\
           0 & -1 & 0 &  0
           \end{pmatrix}
             \begin{pmatrix}
             z_{\mu}z^{\gamma_{1}}\\
             z_{\nu}z^{\gamma_{1}}\\
             z_{\nu}z^{\gamma_{2}}\\
             z_{\mu}z^{\gamma_{2}}
             \end{pmatrix},
        \end{equation}
    where $\gamma_{1}\neq\gamma_{2}$ and $z_{r}^{2}z_{l}^{2}z^{2\gamma}=z_{\mu}z_{\nu}z^{\gamma_{1}}z^{\gamma_{2}}=z^{\beta_{k}}$.
    \end{proposition}

    We need some lemmas.

    If $\beta=(r_{0},r_{1},\ldots,r_{d})$ is a multi-index, then $|\beta|:=r_{0}+r_{1}+\cdots+r_{d}$.

    \begin{lemma}\label{lem4.2}
    Let $\beta=(r_{0},r_{1},\ldots,r_{d})$, $|\beta|=2n_{0}$, be a multi-index and let $\Pi_{\beta}$ be a set of all unordered pairs $\pi=\{z^{\alpha_{i}},z^{\alpha_{j}}\}$ such that $z^{\alpha_{i}}z^{\alpha_{j}}=z^{\beta}$, where $|\alpha_{i}|=|\alpha_{j}|=n_{0}$ and $\deg_{z_{k}} z^{\alpha_{\nu}}\leq n_{k}$, $k=0,1,\ldots,d$. Then $z_{0}^{\delta_{0}}z_{1}^{\delta_{1}}\cdots z_{d}^{\delta_{d}}\in \pi\in \Pi_{\beta}$ if and only if\, $\sum_{k=0}^{d}\delta_{k}=n_{0}$ and
      \begin{equation}\label{eq4.10}
      \max\{(r_{k}-n_{k}),\,0\}\leq \delta_{k}\leq \min\{r_{k},\,n_{k}\}, \quad k=0,1,\ldots,d.
      \end{equation}
    \end{lemma}

    \begin{proof}
    By condition, the inequalities
      $$
      0\leq \delta_{k}\leq \min\{n_{k},\,r_{k}\},\quad  0\leq r_{k}-\delta_{k}\leq \min\{n_{k},\,r_{k}\}
      $$
    holds. From this
      $$
      0\leq \delta_{k}\leq \min\{n_{k},\,r_{k}\},\quad  \max\{(r_{k}-n_{k}),\,0\}\leq \delta_{k}\leq r_{k}.
      $$
    whence (\ref{eq4.10}) follows.
    \end{proof}

    Note that $\pi\in \Pi_{\beta}$ is completely determined by the choice of one of the monomials $z^{\alpha_{i}},z^{\alpha_{j}}$ of the pair $\pi$.

    \begin{definition}\label{def4.1}
    \emph{Let $\Pi_{\beta}$, $\beta=(r_{0},r_{1},\ldots,r_{d})$,  be a set of all unordered pairs $\{z^{\alpha_{\mu}},z^{\alpha_{\nu}}\}$ such that $z^{\alpha_{\mu}}z^{\alpha_{\nu}}=z^{\beta}$, where  $\deg_{z_{k}} z^{\alpha_{\nu}}\leq n_{k}$, $k=0,1,\ldots,d$ and $\deg z^{\alpha_{\nu}}=n_{0}$ for every $\nu$.}

    \noindent
    \emph{If $\pi_{1}=\{z^{\alpha_{i}},\,z^{\alpha_{j}}\}\in \Pi_{\beta}$ and $\deg_{z_{k}}z^{\alpha_{i}}<\min\{r_{k},\,n_{k}\}$, $\deg_{z_{l}} z^{\alpha_{i}}>\max\{(r_{l}-n_{l}),0\}$, then pair}
      \begin{equation}\label{eq4.11}
      \pi_{2}= \left\{\frac{z_{k}}{z_{l}}z^{\alpha_{i}},\,\frac{z_{l}}{z_{k}}z^{\alpha_{j}}\right\}\in \Pi_{\beta}
      \end{equation}
    \emph{is called} an elementary transformation \emph{of pair $\pi_{1}$.
    (\ref{eq4.11}) is uniquely determined by the multi-index} $(0,\ldots,+1_{k},\ldots,-1_{l},\ldots,0)$.
    \end{definition}

    \begin{lemma}\label{lem4.3}
    For any $\pi_{1},\,\pi_{2}\in \Pi_{\beta}$, $\pi_{1}\neq\pi_{2}$, there exists a connecting $\pi_{1}$ and $\pi_{2}$ chain of elements $\pi_{k}\in \Pi_{\beta}$ such that each next element of chain is an elementary transformation of the previous one.
    \end{lemma}

    \begin{proof}
    Suppose $\beta=(r_{0},r_{1},\ldots,r_{d})$. Consider the first monomials of pairs $\pi_{1}$ and $\pi_{2}$. We have $z^{\alpha_{i}}=z_{0}^{\delta_{0}}z_{1}^{\delta_{1}}\cdots z_{d}^{\delta_{d}}\in \pi_{1}$, $z^{\alpha_{j}}=z_{0}^{\nu_{0}}z_{1}^{\nu_{1}}\cdots z_{d}^{\nu_{d}}\in \pi_{2}$ and $z^{\alpha_{i}}\neq z^{\alpha_{j}}$. The monomials $z^{\alpha_{i}}$, $z^{\alpha_{j}}$ satisfy (\ref{eq4.10}).
    If $k_{s}=\nu_{s}-\delta_{s}$, $s=0,1,\ldots,d$, then
      $$
      \sum_{s=0}^{d}k_{s}=0.
      $$
    If $n$ is the sum of positive $k_{s}$, then
    $(k_{0},k_{1},\ldots,k_{d})$ is componentwise sum of $n$ elementary tuples $(m_{0},m_{1},\ldots,m_{d})$, containing only two nonzero components $+1$, $-1$.
    Each $(m_{0},m_{1},\ldots,m_{d})$ corresponds to an elementary transformation  (\ref{eq4.11}). Starting from $\pi_{1}$ at each step we obtain pairs of monomials satisfying (\ref{eq4.10}). At the last step we obtain $\pi_{2}$.
    \end{proof}

    \begin{lemma}\label{lem4.4}
    If\, $\Pi_{\beta}$ contain $m\geq2$ elements, then there exist $m-1$ different twos $\{\pi_{i},\,\pi_{j}\}\in \Pi_{\beta}$ such that the pair $\pi_{j}$ is an elementary transformation of the pair $\pi_{i}$.
    \end{lemma}

    \begin{proof}
    Let us associate the finite graph with the set $\Pi_{\beta}$. The vertices are elements $\pi_{k}\in \Pi_{\beta}$, $k=1,\,\ldots\,,m$. The edges form twos $\{\pi_{i},\,\pi_{j}\}$  connected by an elementary transformation. By Lemma \ref{lem4.3}, the graph is connected. The graph tree contains $m-1$ edges. In the graph tree different edges are incident to different pairs of vertices. Therefore, there are $m-1$ different twos $\{\pi_{i},\,\pi_{j}\}$ connected by an elementary transformation.
    \end{proof}
    \medskip

    \noindent
    \emph{Proof of Proposition \ref{pro4.2}}.
    Taking Remark \ref{rem4.1} into account, we can assume that $\Pi_{\beta}$ contains at least $m\geq2$ elements.
    By Lemma \ref{lem4.4}, there are $m-1$ different twos $\{\pi_{i},\,\pi_{j}\}$ connected by an elementary transformation.
    Let us show that each such twos $\{\pi_{i},\,\pi_{j}\}$ defines a basis matrix of the form (\ref{eq4.8}) or (\ref{eq4.9}). The following cases are possible.
    \smallskip

    \noindent
    (a). $\pi_{1}=\{z^{\alpha_{i}},\,z^{\alpha_{i}}\}$, $\pi_{2}=\{z^{\alpha_{j}},\,z^{\alpha_{k}}\}$.

    \noindent
    Therefore, $z^{\beta}=z^{\alpha_{i}}z^{\alpha_{i}}=z^{\alpha_{j}}z^{\alpha_{k}}$. Since $\pi_{2}$ is obtained by an elementary transformation from $\pi_{1}$, there exist $z_{r}, z_{l}$ such that
      $$
      z^{\alpha_{j}}=\frac{z_{r}}{z_{l}}z^{\alpha_{i}},\quad
      z^{\alpha_{k}}=\frac{z_{l}}{z_{r}}z^{\alpha_{i}}.
      $$
    Let $z^{\gamma}$ be the greatest common divisor of the monomials $z^{\alpha_{i}}$, $z^{\alpha_{j}}$ and $z^{\alpha_{k}}$. Therefore
      $$
      z^{\alpha_{i}}=z_{r}z_{l}z^{\gamma},\quad z^{\alpha_{j}}=z_{r}^{2}z^{\gamma},\quad
      z^{\alpha_{k}}=z_{l}^{2}z^{\gamma}.
      $$
   This triple of monomials defines the basis matrix (\ref{eq4.8}).
   \medskip

   \noindent
   (b). $\pi_{1}=\{z^{\alpha_{i}},\,z^{\alpha_{j}}\}$, $z^{\alpha_{i}}\neq z^{\alpha_{j}}$, $\pi_{2}=\{z^{\alpha_{k}},\,z^{\alpha_{l}}\}$, $z^{\alpha_{k}}\neq z^{\alpha_{l}}$.
   Therefore, $z^{\beta}=z^{\alpha_{i}}z^{\alpha_{j}}=z^{\alpha_{k}}z^{\alpha_{l}}$.
   Since $\pi_{1}$, $\pi_{2}$ be connected by an elementary transformation, then there exist $z_{\mu}, z_{\nu}$ such that
      $$
      z^{\alpha_{k}}=\frac{z_{\mu}}{z_{\nu}}z^{\alpha_{i}},\quad
      z^{\alpha_{l}}=\frac{z_{\nu}}{z_{\mu}}z^{\alpha_{j}}.
      $$
   If $z^{\gamma_{1}}$ is the greatest common divisor of $z^{\alpha_{i}},\,z^{\alpha_{k}}$ and $z^{\gamma_{2}}$ is the greatest common divisor of $z^{\alpha_{j}},\,z^{\alpha_{l}}$, then
     $$
     z^{\alpha_{i}}=z_{\nu}z^{\gamma_{1}},\quad z^{\alpha_{k}}=z_{\mu}z^{\gamma_{1}},\quad
     z^{\alpha_{j}}=z_{\mu}z^{\gamma_{2}},\quad z^{\alpha_{l}}=z_{\nu}z^{\gamma_{2}}.
     $$
   The constructed four monomials define the basis matrix (\ref{eq4.9}). Note that $z^{\gamma_{1}}\neq z^{\gamma_{2}}$. Indeed, if $z^{\gamma_{1}}=z^{\gamma_{2}}$, then $\pi_{1}=\pi_{2}$, which is impossible.

   All twos $\{\pi_{i},\,\pi_{j}\}$ are different. Then corresponding $m-1$ matrices of form (\ref{eq4.8}), (\ref{eq4.9}) are linearly independent. Since $\dim L_{\beta}=m-1$, these matrices form a basis in the subspace $L_{\beta}$.
   \qed
   \medskip

   \noindent
   \emph{Proof of Ambiguity Lemma }(Lemma \ref{lem4.1}).
   Let's put
     $$
     \widetilde{\Psi}(z_{0},z_{1},\ldots,z_{d})=z_{0}^{n_{0}}\Psi(z_{1}/z_{0},\ldots z_{d}/z_{0}).
     $$
   It is easy to see that for $S_{d}$ and $\widetilde{\Psi}(z_{0},z_{1},\ldots,z_{d})$ the relations (\ref{eq4.4}) holds.

   Without loss of generality, we can assume $\widetilde{\Psi}(z)=(z_{d}^{n_{d}}\varphi(\widehat{z}),\,\psi(z))$, where $\widehat{z}=(z_{0},z_{1},\ldots,z_{d-1})$ and $\deg_{z_{d}}\psi(z)=n_{d}-1$. Then
      $$
      S_{d}\frac{\partial^{n_{d}}\widetilde{\Psi}(z)^{T}}{\partial z_{d}^{n_{d}}}=
        \begin{pmatrix}
        S_{11}     & S_{12}\\
        S_{12}^{T} & \widehat{S}_{d}
        \end{pmatrix}
          \begin{pmatrix}
          n_{d}!\varphi(\widehat{z})^{T}\\
          0
          \end{pmatrix}=0.
      $$
   Therefore, $S_{11}=0$, $S_{12}^{T}=0$. From (\ref{eq4.4}) it follows  $\psi(z)\,\widehat{S}_{d}\,\psi(z)^{T}=0$. We rewrite (\ref{eq4.5}) in block form
     \begin{equation}\label{eq4.12}
      \begin{pmatrix}
        S_{11}(\widehat{z})   & S_{12}(\widehat{z}) \\
      S_{12}(\widehat{z})^{T} & S_{22}(\widehat{z})+z_{d}\widehat{S}_{d}
      \end{pmatrix}
        \begin{pmatrix}
        z_{d}^{n_{d}}\varphi(\widehat{z})^{T}\\
        \psi(z)^{T}
        \end{pmatrix}=
          \begin{pmatrix}
          0\\
          0
          \end{pmatrix}.
     \end{equation}
   Relation (\ref{eq4.12}) will be considered as an equation for unknown matrices $S_{ij}(\widehat{z})$, $i,j=1,2$. Since $\psi(z)\,\widehat{S}_{d}\,\psi(z)^{T}=0$, $\deg_{z_{d}}\psi(z)=n_{d}-1$,  then $\widehat{S}_{d}$ is a linear combination of the basis matrices $\widehat{S}_{d,j}$ of form (\ref{eq4.8}), (\ref{eq4.9}). The solution of equation (\ref{eq4.12}) is a linear combination of basic solutions.

   Let us show that the solution of (\ref{eq4.12}) exists for basis matrices $\widehat{S}_{d,j}$.

   For the basis matrix (\ref{eq4.8}), the monomials $z_{r}^{2}z^{\gamma}$, $z_{l}^{2}z^{\gamma}$ and $z_{r}z_{l}z^{\gamma}$ have degree in the variable $z_{d}$ at most $n_{d}-1$. Therefore, there exist monomials $z_{d}z_{r}z^{\gamma},\,z_{d}z_{l}z^{\gamma}\in \widetilde{\Psi}(z)$ such that the nonzero elements of matrix in (\ref{eq4.12}) have the form
     $$
     \left(\begin{array}{cc|ccc}
       0    &   0    &   0    & -z_{l} &   z_{r} \\
       0    &   0    & z_{l}  & -z_{r} &   0     \\
       \hline
       0    & z_{l}  &   0    &    0   &  -z_{d} \\
     -z_{l} & -z_{r} &   0    & 2z_{d} &   0     \\
      z_{r} &   0    & -z_{d} &    0   &   0
     \end{array}\right)
        \left(\begin{array}{c}
        z_{d}z_{r}z^{\gamma}\\
        z_{d}z_{l}z^{\gamma}\\
        \hline
        z_{r}^{2}z^{\gamma} \\
        z_{r}z_{l}z^{\gamma}\\
        z_{l}^{2}z^{\gamma}
        \end{array}\right)=0.
     $$
   Similarly, for basis matrix of form (\ref{eq4.9}) and monomials $z_{d}z^{\gamma_{2}},\,z_{d}z^{\gamma_{1}}\in \Psi(z)$ we have
    $$
    \left(\begin{array}{cc|cccc}
        0    &     0    & -z_{\nu} &  z_{\mu} &     0    &    0   \\
        0    &     0    &     0    &      0   & -z_{\mu} & z_{\nu}\\
        \hline
    -z_{\nu} &     0    &     0    &      0   &  z_{d}   &    0   \\
     z_{\mu} &     0    &     0    &      0   &     0    & -z_{d} \\
        0    & -z_{\mu} &  z_{d}   &      0   &     0    &    0   \\
        0    &  z_{\nu} &     0    & -z_{d}   &     0    &    0
    \end{array}\right)
      \left(\begin{array}{c}
        z_{d}z^{\gamma_{2}}\\
        z_{d}z^{\gamma_{1}}\\
        \hline
      z_{\mu}z^{\gamma_{1}}\\
      z_{\nu}z^{\gamma_{1}}\\
      z_{\nu}z^{\gamma_{2}}\\
      z_{\mu}z^{\gamma_{2}}
      \end{array}\right)=0.
    $$
  Therefore, for (\ref{eq4.12}) we obtain
     $$
     (z_{0}S_{0}+z_{1}S_{1}+\cdots+z_{d-1}S_{d-1}+z_{d}S_{d})\widetilde{\Psi}^{T}(z_{0},z_{1},
     \ldots,z_{d})=0.
     $$
  Assuming $z_{0}=1$, we get (\ref{eq4.5}).
  \qed
  \medskip

  \noindent
  \emph{Proof of assertion \emph{(ii)} of Theorem \ref{the3.1}}. By (i) of Theorem \ref{the3.1}, we have
     $$
     q(\zeta)p(z)=\Psi(\zeta)(B_{0}+z_{1}B_{1}+\cdots+z_{d}B_{d})\Psi(z)^{T},
     $$
  and
     $$
     F_{d}(z)=q(z)\frac{\partial p(z)}{\partial z_{d}}-p(z)\frac{\partial q(z)}{\partial z_{d}}=
     \Psi(z)B_{d}\Psi(z)^{T}.
     $$
  Suppose that  $F_{d}(z)$ is a SOS polynomial.
  By Theorem 1 (\cite{uj21}), there exists a positive semidefinite real $N\times N$ matrix $A_{d}$ such that
     $$
     F_{d}(z)=q(z)\frac{\partial p(z)}{\partial z_{d}}-p(z)\frac{\partial q(z)}{\partial z_{d}}=
     \Psi(z)A_{d}\Psi(z)^{T}.
     $$

   Since $F_{d}(x)\geq 0$ for $x\in \mathbb{R}^{d}$, then $\deg_{z_{d}} F_{d}(z)\leq2(n_{d}-1)$. If $\deg_{z_{d}}a_{ii}z^{\alpha_{i}}z^{\alpha_{i}}=2n_{d}$, then $a_{ii}=0$. From $A_{d}\geq0$ we get
     $$
     A_{d}\frac{\partial^{n_{d}}\Psi(z)^{T}}{\partial z_{d}^{n_{d}}}=0.
     $$
   Therefore, real symmetric matrix $S_{d}=A_{d}-B_{d}$ satisfies the conditions of  Ambiguity Lemma. Then there exist real symmetric matrices $S_{0},\,\ldots\,,S_{d-1}$ for which
     $$
     \Psi(\zeta)(S_{0}+z_{1}S_{1}+\cdots+z_{d-1}S_{d-1}+z_{d}S_{d})\Psi(z)^{T}=0.
     $$
  For the matrix pencil $A(z)=B(z)+S(z)$ we get
     $$
     q(\zeta)p(z)=\Psi(\zeta)(A_{0}+z_{1}A_{1}+\cdots+z_{d}A_{d})\Psi(z)^{T},
     $$
  where $A_{d}=B_{d}+S_{d}\geq 0$.
  \qed

  \begin{corollary}\label{cor4.1}
  Let $f(z)=p(z)/q(z)$, $z\in \mathbb{C}^{d}$, be a rational real function. If
      $$
      F_{d}(z)=q(z)\frac{\partial p(z)}{\partial z_{d}}-p(z)\frac{\partial q(z)}{\partial z_{d}}
      $$
  is a PSD polynomial, then
      $$
      f(z)=\left[\pi(A_{0}+z_{1}A_{1}+\cdots+z_{d}A_{d})^{-1}\pi^{T}\right]^{-1},
      $$
  where $\pi=(1,0,\ldots,0)$, $A_{j}$, $j=0,1,\ldots,d$, are real symmetric $M\times M$-matrices, and $A_{d}\geq 0$.
  \end{corollary}

  \begin{proof}
  By Artin's theorem, there exists polynomial\; $s(z)$\; such that\; $s(z)^{2}F_{d}(z)$ is SOS.
  For $q(z)s(z)$, $p(z)s(z)$, by Theorem \ref{the3.1} (ii) we obtain
     $$
     q(\zeta)s(\zeta)p(z)s(z)=\Psi(\zeta)(B_{0}+z_{1}B_{1}+\cdots+z_{d}B_{d})\Psi(z)^{T},
     $$
  where $B_{d}\geq 0$. Now it suffices to apply Corollary \ref{cor3.1} and Remark \ref{rem3.1}.
  \end{proof}

  \section{Main Theorem}\label{s:5}

  Let $z\in\mathbb{C}^{d}$. Denote by $\mathfrak{N}_{d}$ the set of rational functions $f(z)$ with real coefficients, holomorphic and satisfying the condition
     \begin{equation}\label{eq5.1}
     \MYim f(z)\geq 0
     \end{equation}
  in region $\MYim z_{1}>0,\ldots,\MYim z_{d}>0$.

  \begin{lemma}\label{lem5.1}
  If $f(z)\in\mathfrak{N}_{d}$, then
     $$
     F_{j}(z)=q(z)\frac{\partial p}{\partial z_{j}}(z)-p(z)\frac{\partial q}{\partial z_{j}}(z),
     \quad j=1,\ldots,d
     $$
  are PSD polynomials.
  \end{lemma}

  \begin{proof}
  Suppose $j=1$. Let us $\varphi(\zeta)=f(\zeta,\widehat{x})$. If $\widehat{x}=(x_{2},\ldots,x_{d})\in\mathbb{R}^{d-1}$, then
          $\mbox{Im}\,\varphi (\zeta)\geq 0\enskip \mbox{for}\enskip \mbox{Im}\,\zeta>0$
  and $\mbox{Im}\,\varphi (\zeta)=0\enskip \mbox{for}\enskip\mbox{Im}\,\zeta=0$.
  Therefore, $d\varphi(\zeta)/d\zeta \mid_{\zeta\in\mathbb{R}}\geq 0$. From this
      \medskip

      \qquad \qquad $F_{1}(x)=q(x)^{2}d\varphi(\zeta)/ d\zeta|_{\zeta=x_{1}}\geq0, \quad x\in\mathbb{R}^{d}$.
   \end{proof}

   The following lemma is a modification of Theorem 1.3.2 (\cite{uj22}).

  \begin{lemma}\label{lem5.2}
   Let $h(z),s(z)\in \mathbb{R}[z_{1},\ldots,z_{d}]$ $(d>1)$ be coprime polynomials, $\partial s/\partial z_{d}\neq 0$ and $\mathbb{K}=\mathbb{R}$ or $\mathbb{C}$. Let, further, $\widehat{z}\in\mathbb{K}^{d-1}$, $z_{d}\in\mathbb{C}$ and $V(\widehat{z})\subseteq\mathbb{K}^{d-1}$, $U(z_{d})\subseteq\mathbb{C}$ be arbitrary neighborhoods of the points $\widehat{z}$, $z_{d}$.

   If $h(\widehat{z},z_{d})=s(\widehat{z},z_{d})=0$ then there exist $\widehat{z}'\in V(\widehat{z})$, $z'_{d}\in U(z_{d})$ such that $s(\widehat{z}', z'_{d})=0$ and $h(\widehat{z}', z'_{d})\neq 0$.
   \end{lemma}

   \begin{proof}
   Note that there is a neighborhood $\Omega(\widehat{z})\subseteq V(\widehat{z})$  such that for each $\widehat{z}'\in \Omega(\widehat{z})$ there is $z'_{d}\in U(z_{d})$ for which $s(\widehat{z}',z'_{d})=0$.

   Let $F=\mathbb{R}(z_{1},\ldots,z_{d-1})$ be the field of quotients  of the ring $R=\mathbb{R}[z_{1},\ldots,z_{d-1}]$. Polynomials $h(z),s(z)$ will be considered as elements of the ring $R[z_{d}]\subset F[z_{d}]$.
   Since $h(z),s(z)$ are coprime in $\mathbb{R}[z_{1},\ldots,z_{d}]=R[z_{d}]$, then, by the Gauss lemma, $h(z),s(z)$ are coprime in $F[z_{d}]$. There are elements $\widehat{a},\widehat{b}\in F[z_{d}]$ for which $\widehat{a}h+\widehat{b}s=1$. We multiply this identity by a nonzero polynomial $c(\widehat{z})\in R=\mathbb{R}[z_{1},\ldots,z_{d-1}]$ such that $\widehat{a}c=a\in R[z_{d}]$, $\widehat{b}c=b\in R[z_{d}]$. Then
     $$
     a(\widehat{z},z_{d})h(\widehat{z},z_{d})+b(\widehat{z},z_{d})s(\widehat{z},z_{d})=c(\widehat{z}).
     $$
   Since $c(\widehat{z})$ is a nonzero polynomial, then there is $\widehat{z}'\in \Omega(\widehat{z})\subseteq V(\widehat{z})$ such that $c(\widehat{z}')\neq0$. There exists $z'_{d}\in U(z_{d})$ for which $s(\widehat{z}',z'_{d})=0$. From the last identity follows $h(\widehat{z}',z'_{d})\neq0$.
   \end{proof}

  \begin{theorem}\label{the5.1} \emph{\textbf{(Main Theorem)}}.
  If $f(z)=p(z)/q(z)\in \mathfrak{N}_{d}$, then
    \begin{equation}\label{eq5.2}
    W_{j}[q(z),p(z)]=q(z)\frac{\partial p(z)}{\partial z_{j}}-p(z)\frac{\partial q(z)}{\partial z_{j}},\quad j=1,\ldots,d
    \end{equation}
  are sums of squares of polynomials.
  \end{theorem}

  \begin{proof}
  It can be considered that $f(z)\neq0$.
  By Lemma \ref{lem5.1}, $W_{j}[q,\,p]$ are PSD polynomials.

  Suppose $j=d$. Let $s(z)$ be Artin's minimum denominator for $F_{d}(z)$. By Theorem \ref{the3.1} (ii), there are real symmetric matrices $A_{0},\ldots,A_{d-1},A_{d}$, where $A_{d}\geq0$, such that
     $$
     s(\zeta)q(\zeta)s(z)p(z)=\Psi(\zeta)(A_{0}+z_{1}A_{1}+\cdots+z_{d}A_{d})\Psi(z)^{T}.
     $$
  This implies
     \begin{equation}\label{eq5.4}
     f(z)=\frac{\Psi(\zeta)}{s(\zeta)q(\zeta)}(A_{0}+z_{1}A_{1}+\cdots+z_{d}A_{d})
     \frac{\Psi(z)^{T}}{s(z)p(z)},
     \end{equation}
       \begin{equation}\label{eq5.5}
       s^{2}(z)F_{d}(z)=\Psi(z)A_{d}\Psi(z)^{T}=H(z)H(z)^{T},
       \end{equation}
  where $H(z)=(h_{1}(z),\ldots,h_{k}(z))$, $h_{i}\in \mathbb{R}[z_{1},\ldots,z_{d}]$.
  Let $s=s_{1}\cdots s_{r}$ be a decomposition of $s(z)$ into irreducible factors.
  Each irreducible factor $s_{j}(z)$ cannot be a divisor of all elements $h_{i}(z)$, otherwise $s(z)/s_{j}(z)$ is also Artin's denominator of $F_{d}$. This contradicts the minimality of $s(z)$. Without loss of generality, we can assume $h_{1}(z)$ and $s_{1}(z)$ are coprime polynomials.

  Assuming $\zeta=\overline{z}=(\overline{z}_{1},\ldots,\overline{z}_{d})$, from (\ref{eq5.4}) we get
    \begin{equation}\label{eq5.6}
    \MYim f(z)=\sum_{j=1}^{d-1}\MYim z_{j}\frac{\Psi(z)A_{j}\Psi(z)^{\ast}}{|q(z)|^{2}|s(z)|^{2}}+
    \MYim z_{d}\frac{H(z)H(z)^{\ast}}{|q(z)|^{2}|s(z)|^{2}},
    \end{equation}
  where $H(z)^{\ast}:=\overline{H(z)}^{T}$.

  Suppose $s(z)$ is different from a constant. The following cases are possible.
  \medskip

  \noindent
  (a). $\partial F_{d}(z)/\partial z_{d}\neq 0$. By Proposition \ref{pro2.3}, $s(x',z_{d})=0$ for $x'\in\mathbb{R}^{d-1}$, $\MYim z_{d}>0$. If $h_{1}(x',z_{d})=0$, then
  by Lemma \ref{lem5.2}, there exists $z'=(x,\zeta)$, $\MYim \zeta>0$, $x\in \mathbb{R}^{d-1}$ such that $h_{1}(z')\neq 0$, $s_{1}(z')=0$. From (\ref{eq5.6}) we get
     $$
     \MYim f(x,z_{d})=\MYim z_{d}\frac{H(x,z_{d})H(x,z_{d})^{\ast}}{|q(x,z_{d})|^{2}|s(x,z_{d})|^{2}},\quad x\in \mathbb{R}^{d-1}.
     $$
  This implies $\lim_{z_{d}\rightarrow z'_{d}}\MYim f(x,z_{d})=\infty$. Since $f(x,z_{d})$ is holomorphic for $\MYim z_{d}>0$, this is impossible.
  \medskip

  \noindent
  (b). $\partial F_{d}(z)/\partial z_{d}=0$. In this case, $H(z)$ and $s(z)$ do not depend on the variable $z_{d}$. By Corollary \ref{cor2.1} and Lemma \ref{lem5.2}, there exists $\widehat{z}=(\zeta_{1},\ldots,\zeta_{d-1})$, $\MYim \zeta_{j}>0$, $j=1,\ldots,d-1$, for which
  $h_{1}(\widehat{z})\neq0$, $s_{1}(\widehat{z})=0$.

   For fixed $z_{d}=x_{d}\in \mathbb{R}$ and $\widetilde{z}$ in some neighborhood $\Omega$ (lying in the upper poly-half-plane) of the point $\widehat{z}$ from (\ref{eq5.1}) and (\ref{eq5.6}) we get
    $$
    \MYim f(\widetilde{z},x_{d})=\sum_{j=1}^{d-1}\MYim z_{j}\frac{\Psi(\widetilde{z},x_{d})A_{j}\Psi(\widetilde{z},x_{d})^{\ast}}
    {|q(\widetilde{z},x_{d})|^{2}|s(\widetilde{z})|^{2}}>0
    $$
   For sufficiently small fixed $y_{d}>0$ and $\widetilde{z}$ in some neighborhood $\Omega'\subseteq \Omega$ this implies
    \begin{equation}\label{eq5.7}
    \sum_{j=1}^{d-1}\MYim z_{j}\frac{\Psi(\widetilde{z},x_{d}+iy_{d})A_{j}\Psi(\widetilde{z},x_{d}+iy_{d})^{\ast}}
    {|q(\widetilde{z},x_{d}+iy_{d})|^{2}|s(\widetilde{z})|^{2}}>0
    \end{equation}
   Then from (\ref{eq5.6}) we obtain
    \begin{multline}\nonumber
    \quad\quad \MYim f(\widetilde{z},x_{d}+iy_{d})=\\
    \sum_{j=1}^{d-1}\MYim z_{j}\frac{\Psi(\widetilde{z},x_{d}+iy_{d})A_{j}\Psi(\widetilde{z},x_{d}+iy_{d})^{\ast}}
    {|q(\widetilde{z},x_{d}+iy_{d})|^{2}|s(\widetilde{z})|^{2}}+
       y_{d}\frac{H(\widetilde{z})H(\widetilde{z})^{\ast}}
       {|q(x,x_{d}+iy_{d})|^{2}|s(\widetilde{z})|^{2}}.
    \end{multline}
  The first term is positive as $\widetilde{z}\rightarrow \widehat{z}$, and the second increases indefinitely. Since $f(z)$ is holomorphic in region $\MYim z_{1}>0,\ldots,\MYim z_{d}>0$, this is impossible.

  The assumption $s(z)$ is different from a constant leads to a contradiction. Then Artin's minimal denominator $s(z)\equiv const$ and $F_{d}(z)$ is the sum of squares of polynomials.
  \end{proof}

\end{document}